\documentclass[12pt]{article}

\usepackage{amssymb,
amsmath, amscd
}

\usepackage{fullpage}
\usepackage{titling}

\usepackage[matrix,arrow,curve,frame,arc]{xy}
\usepackage{tikz}
\usepackage{color}

\newcommand{\La}{\Lambda}

\newcommand{\cQ}{\mathcal{Q}}
\newcommand{\cB}{\mathcal{B}}
\newcommand{\cD}{\mathcal{D}}
\newcommand{\vf}{\varphi}
\newcommand{\C}{\mathbb{C}}
\newcommand{\bZ}{\mathbb{Z}}
\newcommand{\N}{\mathbb{N}}
\newcommand{\cT}{\mathcal{T}}
\newtheorem{defi}{Definition}[section]
\newtheorem{rem}[defi]{Remark}
\newtheorem{prop}[defi]{Proposition}

\newtheorem{lemma}[defi]{Lemma}
\newtheorem{cor}[defi]{Corollary}
\newtheorem{thm}[defi]{Theorem}
\newtheorem{expl}[defi]{Example}

\newcommand{\cX}{\mathcal{X}}

\newcommand{\ul}{\underline}
\newcommand{\ve}{\varepsilon}
\newcommand{\wt}{\widetilde}

\DeclareRobustCommand{\authorthing}{
Karin Erdmann\\
Mathematical Institute, ROQ, Oxford OX2 6GG, UK\\
{\tt erdmann@maths.ox.ac.uk}
}
\date{}
\author{\authorthing}

\title{On $\ell$-parabolic Hecke algebras of symmetric groups}

\begin{document}

\maketitle


\abstract
\footnote{2010 Mathematics Subject classification 20C08, 16E40, 16G10, 16G70,
16P10, 16P20.}
{Let $H = H_q(n)$ be the Hecke algebra of the symmetric group
of degree $n$,  over a field of arbitrary characteristic, and 
where $q$ is a primitive $\ell$-th root of unity in $K$.
Let $H_{\rho}$ be an $\ell$-parabolic subalgebra
of $H$. 
We give an elementary explicit construction for the basic algebra
of a non-simple block of $H_{\rho}$. 
We  also discuss homological properties of $H_{\rho}$-modules, 
in particular existence of varieties for modules, and some consequences.}

\section{Introduction}

Let $H=H_q(n)$ be the Hecke algebra of a symmetric group over some field $K$, where $q\in K$ is a primitive
$\ell$-th root of unity with $\ell > 1$, and where $K$ has characteristic $p\geq 0$. One would like to understand
homological properties of $H$, and how they may differ from those of group algebras.

Standard parabolic subalgebras of $H$ are an analog of subgroup algebras for group algebras. In \cite{DD} a vertex theory was
developed; and it is now known that a vertex of an indecomposable $H$-module is $\ell-p$ parabolic (see \cite{DD}, \cite{D}, and \cite{W} for the general case). 
For any composition $\rho$ of $n$, there is an associated  parabolic subalgebra $H_{\rho}$ of $H$, which is isomorphic to  the Hecke algebra of the
standard Young subgroup $W_{\rho}$ of the symmetric group  of degree $n$. The algebra $H_{\rho}$ is $\ell$-parabolic if all parts of $\rho$ are $\ell$ or $1$, and
it is $\ell- p$ parabolic if all parts are $1, \ell, \ell p, \ell p^2, \ldots $.
This suggests that $\ell-p$ parabolic subalgebras should be the analogs of group algebras of $p$-groups, and that  $\ell$-parabolic subalgebras
should play a  role similar to that of group algebras of elementary abelian $p$-groups.
\medskip

Here we study $\ell$-parabolic subalgebras $H_{\rho}$ of $H$, the smallest one
is $H_q(\ell)$.  
We give an elementary and explicit
construction for the basic algebra of the principal block of $H_q(\ell)$, the unique non-simple block. (The structure is known, 
over characteristic $0$ it may be found in \cite{Ge} or \cite{U}
by a rather indirect approach.)
By taking
  tensor powers, this gives explicit basic algebras for an arbitrary  non-simple block of any $\ell$-parabolic subalgebra $H_{\rho}$.
Recall that, in general,  a basic algebra $A$ of an algebra $\La$ 
is the smallest algebra which is Morita equivalent to $\La$, it is unique up to isomorphism. When the field is large enough, the simple $A$-modules are one-dimensional. For many problems working, with
a basic algebra has great advantages.

In the second part of the paper, we study homological
properties of $\ell$-parabolic Hecke algebras. 
Given a module $M$ of some algebra $\La$, 
it is a basic question  whether or not it is projective. 
More generally, one wants to know the complexity of $M$, that is the rate
of growth of a minimal projective resolution.  To answer these, and
understand other  invariants, one can exploit cohomological
support varieties of modules, or perhaps rank varieties, when they exist.

We show that $\ell$-parabolic Hecke algebras have  a theory
of support varieties
constructed via Hochschild cohomology; this can be mostly extracted from
existing literature.
As a consequence, one also can describe
possible tree classes of Auslander-Reiten components.

Rank varieties
were originally introduced by J. Carlson \cite{C} for
group algebras of elementary abelian $p$-groups over fields
of characteristic $p$.
They were generalized to a class of quantum complete intersections in \cite{BEH} (see also \cite{Sch} for a more detailed account). 
We show that any non-simple block of an $\ell$-parabolic Hecke algebra
has a subalgebra $R$ which controls directly projectivity of modules in the block. 
By \cite{BEH}, this algebra $R$ has a theory of rank varieties; though
it does not  seem to extend to $H_{\rho}$ in general.

\bigskip
Let $\cX$ be the category of indecomposable $H$-modules which are not projective restricted to $H_{\rho}$, it is a union of stable Auslander-Reiten components. We show that one can also describe the possible tree classes
of components containing modules in $\cX$. 

Section 2 contains background on  Hecke algebras and basic algebras. Section 3 constructs the basic algebra of $H_q(\ell)$, and in Section 4 we discuss
homological properties. Analogues for group algebras of symmetric groups
also hold,  with $p$ instead of $\ell$.

Our approach to  Hecke algebras  is based on \cite{DJ}, and 
for background on homological properties we refer to 
\cite{Be}. 
We work with right modules.

\medskip

{\bf Acknowledgement.} \  Most of 
this material is based on work supported by the
National Science Foundation under Grant No. DMS-1440140 while the author was in
residence at the Mathematical Sciences Research Institute in Berkely, California, during part of the Spring 2018 semester. The author thanks Dave Benson and Dan Nakano  for discussions related to this material, and thanks to the referee.

\bigskip

\parindent0pt

\section{Preliminaries}

{\bf 2.1 } \ Throughout, let 
$n$ be a positive integer, and let $W$ be the symmetric group
on $\{ 1, 2, \ldots, n\}$. The following is based on \cite{DJ}.

Let $K$ be a field, and let $q$ be a primitive $\ell$-th root of unity in 
$K$ where $\ell \geq 2$. The Iwahori Hecke algebra 
$H= H_q(n)$ has generators $T_1, \ldots, T_{n-1}$ and relations\\
$$T_iT_{i+1}T_i = T_{i+1}T_iT_{i+1} \ \ (1\leq i < n-1), \ \ \ T_iT_j= T_jT_i
\ \ (|i-j|\geq 2)$$
$$T_i^2 = (q-1)T_i + qI \ \ (1\leq i \leq n - 1).
$$
It has basis $\{ T_w: w\in W\}$. If $w= (i \ i+1)=:s_i$ then 
$T_w= T_{i}$. The algebra $H_q(n)$ is a deformation of the group algebra  $KW$. 
It is a symmetric algebra but is not a Hopf algebra.

Let $\cB$ be the set of basic transpositions, that is $\cB=\{ s_1, \ldots, s_{n-1}\}$. Recall that the length of an element $w\in W$ is the minimal $k\geq 0$ such that
$w$ can be written as $w= v_1v_2\ldots v_k$ for $v_i\in \cB$. If so then
$T_w = T_{v_1}\ldots T_{v_k}$. 

\medskip

A composition of $n$ is a tuple of non-negative integers
$\lambda = (\lambda_1, \lambda_2, \ldots, \lambda_a)$ such that $\sum_i \lambda_i=n$.
For such a composition, let  $W_{\lambda}$ be the corresponding standard
Young subgroup of $W$ , it is the direct product
$W_{[1, \lambda_1]}\times W_{[\lambda_1+1, \lambda_1+\lambda_2]} \times \ldots
\times W_{[\sum_{i=1}^{a-1}\lambda_i+1, n]}$.
Here $[u, v] = \{ x\in \N\mid u\leq x\leq v\}$, we refer to such a set as 
an {\it interval}. 
The parabolic subalgebra $H_{\lambda}$ is the
Hecke algebra of the Young subgroup $W_{\lambda}$, it is 
isomorphic to the tensor product of
 algebras $H_q(\lambda_i)$ with disjoint supports.

\medskip

{\bf 2.2} \ 
(a) \ In order to work with a parabolic subalgebra $H_{\lambda}$ of $H$, 
one needs to use distinguished coset representatives. That is, there is
a set $\cD_{\lambda}$ of permutations such that 
$W=\bigcup_{d\in \cD_{\lambda}} W_{\lambda}d$  (the disjoint union),
and moreover $d$ is the unique element of minimal length in its coset.
If so, then  
 $(\cD_{\lambda})^{-1} = \{ d^{-1}\mid d\in \cD_{\lambda}\}$ 
 is a system of left coset representatives, where each element 
 is the unique element of minimal length in its coset.

To compute 
the elements in $\cD_{\lambda}$, let
$t^{\lambda}$ be
the standard $\lambda$-tableau in which the numbers  $1, 2, \ldots, n$ appear
in order along  successive rows.
With this, $\cD_{\lambda}$ is the set of permutations $d$ such that
the tableau $t^{\lambda}d$ is a row standard tableau.
\medskip

(b) \ In order to work with pairs of parabolic subalgebras $H_{\lambda}$ and
$H_{\mu}$ for compositions $\lambda, \mu$, we need double cosets. 
Let $\cD_{\lambda,\mu} = \cD_{\lambda}\cap \cD_{\mu}^{-1}$, this is
a system of distinguished  
$W_{\lambda}-W_{\mu}$ double coset representatives in $W$, and 
if $d\in \cD_{\lambda, \mu}$ then $d$ is the unique
element of minimal length in the double coset $W_{\lambda}dW_{\mu}$.
In 1.7 of \cite{DJ} it is described how to compute $\cD_{\lambda, \mu}$.
We will use a modification below.

\medskip

{\bf 2.3 } 
\ Let $R$ be a non-empty interval contained in $\{ 1, 2, \ldots, n \}$,
suppose $W_R$ is the group of permutations of $R$, and   $H_q(R)$ is the Hecke algebra of $W_R$.
In  $H_q(R)$
 we consider 
$$x_R:= \sum_{w\in W_R} T_w, \ \ y_R:= \sum_{w\in W_R} \ (-q)^{-l(w)}T_w.
$$
Then  $x_R$ and $y_R$ span the trivial, and the sign module of $H_q(R)$, 
that is 
$$x_{R}T_w= q^{l(w)}x_{R}  \ \mbox{ and} \  y_RT_w = (-1)^{l(w)}y_R \ \ (w\in W_R).
$$
The elements $x_R, y_R$ are central in $H_q(R)$, and
their product is zero if $R$ has size $\geq 2$ 
and $\ell \geq 3$.

We have
$$x_{R}^2 = (\sum_{w\in W_R} q^{l(w)})x_R \ \ \mbox{ and } \ \ 
y_R^2 = ((\sum_{w\in W_R} q^{-1})^{l(w)})y_R.
$$
As is well-known,
if $u$ is a variable then
$\sum_{w\in W_R} u^{l(w)} = \prod_{m=1}^r [m]_u =: [r]_u!
$
where $[m]_u  =  1 + u + \ldots + u^{m-1}=\frac{1-u^m}{1-u}$.
Assume now that $r < \ell$, then if we substitute $u=q$ or $q^{-1}$  we
 get a non-zero element in $K$. This means that
we have idempotents
$c_Rx_R$ and $\hat{c}_Ry_R$
where $c_R= \frac{1}{[r]_q!}$ and $\hat{c}_R=\frac{1}{[r]_{\hat{q}}!}$, setting
$\hat{q}=q^{-1}$. 
Note that the condition $r<\ell$ is essential; if $r=\ell$ then $x_R^2=0=y_R^2$.

 \medskip

More generally, if $\rho$ is a composition of $n$ , define $x_{\rho}$, or $y_{\rho}$, 
to be the product of the $x_{\rho_i}$, or $y_{\rho_i}$, 
where the $\rho_i$ are the support sets of the
standard Young subgroup $W_{\rho}$.  Then we have
$$x_{\rho}^2  = c_{\rho} x_{\rho}, \ y_{\rho}^2 = \hat{c}_{\rho}y_{\rho}
$$
where $c_{\rho}$ and $\hat{c}_{\rho}$ are
non-zero scalars.
Suppose $\rho$ is a composition  of $R$ and $|R|<\ell$, then 
$x_{\rho}x_R$ is a non-zero scalar multiple of $x_R$, and $y_{\rho}y_R$
is a non-zero scalar multiple of $y_R$. If $\rho_i\geq 2$ for at least one
$i$ then $y_{\rho}x_R=0=x_Ry_{\rho}$ and $x_{\rho}y_R=0=y_Rx_{\rho}$.

\medskip

{\bf 2.4 }
\  Let $\#$ be the automorphism of $H$ given by
$T_w\mapsto (-q)^{l(w)}(T_{w^{-1}})^{-1}$. Twisting  by this automorphism 
interchanges  $x_R$ and $y_R$ for 
any interval.

There is also an anti-involution $*$ on $H$ defined on the
basis, by  $T_w \mapsto T_{w^{-1}}$, it fixes $x_R$ and $y_R$ for any interval
$R$.  
 Recall from \cite{DJ}
the symmetric bilinear form 
$$(T_u, T_v) = \left\{\begin{array}{cc} q^{l(u)},  & u=v\cr 0, & \mbox{ else.}
\end{array}
\right.
$$
This satisfies
$$(h_1h_2, h_3) = (h_1, h_3\cdot h_2^*) \ \ \mbox{\ and} \ 
(h_1, h_2h_3) = (h_2^*h_1, h_3) \ (h_i\in H).$$
If one  defines
$f(h_1, h_2):= (h_1, h_2^*)
$
then, as one can check, 
$f$ is a symmetrizing form. That is, it is symmetric, associative,
and non-degenerate, and it shows that $H$ is a symmetric algebra.

\medskip

{\bf 2.5  } \  The combinatorics for representations of Hecke algebras
is similar to that  for symmetric groups, with $\ell$ instead of $p$, as it is proved in \cite{DJ} and \cite{DJ2}, see also \cite{Ma}.
For each partition $\lambda$ of $n$, there is a Specht module $S^{\lambda}$. When $\lambda$ is $\ell$-regular,
this has a unique simple quotient $D^{\lambda}$, and these are the simple
$H$-modules,  up to isomorphism. 
The blocks are labelled by $\ell$-cores $\kappa$ and weight $w$, for $n= |\kappa| + w\ell$. 
The principal block is the block  containing the trivial module, that is
$S^{\lambda}$ for $\lambda = (n)$. 
For our construction we  only use that when $n=\ell$, the principal block
has precisely $\ell-1$ simple modules.

\medskip

{\bf 2.6 } We recall facts about basic algebras, for details see for example
\cite{ASS}.  Assume $A$ is any finite-dimensional $K$-algebra, a basic algebra associated to $A$ is an algebra
$eAe$ where $e$ is an idempotent of the form $e=e_1 + \ldots + e_m$ with $e_i$ orthogonal primitive idempotents, such
that $e_iA$ is not isomorphic to $e_jA$ for $i\neq j$, and such that every indecomposable projective $A$-module is isomorphic
to some $e_iA$. 
If so, for large enough  $K$, 
the $eAe$ is isomorphic to an algebra $KQ/I$ where $Q$ is a 
unique quiver (directed graph),  $KQ$ is the path algebra of $Q$, and 
$I$ is an admissible ideal of $KQ$ (that is $KQ^N\leq I \leq KQ^2$ where $KQ^r$ is the span of the paths of lengths $\geq r$. 

The vertices of $Q$ correspond to the idempotents $e_i$, and the arrows $i\to j$ are in bijection with a basis of the vector space
$e_i({\rm rad} A/{\rm rad}^2A)e_j$.  
The connected components of $Q$ are in bijection with the blocks of the algebra $A$ (ie the indecomposable direct summands of $A$
as an algebra). For a connected component  of $Q$ with vertices $V \subseteq \{ 1, 2, \ldots, m\}$ set $e_V:= \sum_{i\in V} e_i$, then
$e_VAe_V$ is a basic algebra for a block of $A$.

\vspace*{0.5cm}

\section{The principal block of $H_q(\ell)$}

We take $H:= H_q(\ell)$. 
The principal block of this algebra 
is the block containing the trivial module. 

Suppose  $\ell = 2$, that is $q=-1$. Then we have $(T_1+1)^2=0$ and the algebra map 
$K[X]\to H$ sending $X\mapsto T_1+1$ induces an isomorphism 
$K[X]/(X^2) \cong H$. The algebra  $H$ is its own basic algebra.
We assume from now  that $\ell \geq 3$.

It is known (see \cite{Ge} or \cite{U} for the case $K=\C$) that a basic algebra of this block  is a Brauer tree algebra,   isomorphic to $K\cQ/I$ where
$\cQ$ is of the form

$$\xymatrix{
      1  \ar@<.5ex>[r]  &
      2 \ar@<.5ex>[l]  \ar@<.5ex>[r]  &
      \cdots \ar@<.5ex>[l] \ar@<.5ex>[r] &
      l-2 \ar@<.5ex>[l] \ar@<.5ex>[r] &
      l-1 \ar@<.5ex>[l],  &} 
$$
and if we label the arrow $i \to i+1$ by $\alpha_i$, and the arrow $i+1\to i$ by $\beta_i$, 
then the ideal
$I$ is generated by  the relations 
$$\beta_i\alpha_i = \alpha_{i+1}\beta_{i+1} \ \ (1\leq i \leq \ell-2), \ \ 
\alpha_1\beta_1\alpha_1, \  \beta_{\ell-2}\alpha_{\ell - 2}\beta_{\ell - 2}.
\leqno{(3.1)}$$
This is a Brauer tree algebra. 
In order to give a 
direct construction of the  basic algebra, as a subalgebra of $H$, we
start by finding a set of 
$\ell-1$ orthogonal idempotents.
They will be of the form 
$x_R y_S$ where $R, S$ are non-empty intervals
which form a partition of $\{ 1, 2, \ldots, \ell \}$. 
As explained in 2.3, we have idempotents  which
are non-zero scalar multiples of $x_R, y_S$, and since
$x_R, y_S$ commute, then also  $x_Ry_S$ is a non-zero scalar multiple of 
an idempotent in $H$. 

\bigskip

\subsection{Idempotents}

Such idempotents have certain factorisations.

Let $R$ be an interval with  $1 \leq |R| < \ell$, and let 
 $W= W_R$. Let $\rho$ be a composition of $R$ and
$$W = \bigcup_{d\in \cD_{\rho}} W_{\rho}d. 
$$

\begin{lemma}\label{lem:3.1} With this, we have\\
(a)  $y_R = y_{\rho}\sigma = \sigma^*y_{\rho}$
where $\sigma = \sum_{d\in \cD_{\rho}} (-q)^{-l(d)}T_d$. \\
(b) \ $x_R = x_{\rho}\psi = \psi^*x_{\rho}$ where
$\psi = \sum_{d\in \cD_{\rho}} T_d$.
\end{lemma}

\bigskip

{\it Proof } (a) \ Let $w\in W_R$, then $w=hd_i$ for $h\in W_{\rho}$ and 
$d_i\in \cD_{\rho}$, and then $l(w) = l(h) + l(d_i)$ and 
$T_w = T_hT_{d_i}$. Therefore 
we can write $y_R$ as 
$$\sum_{w} (-q)^{-l(h)-l(d_i)}T_hT_{d_i}
$$ and this gives the required factorisation. 
We apply the anti-involution $*$ (see 2.4) and get 
$y_R= y_R^* = \sigma^*y_{\rho}$. 
Part (b) is similar.
$\Box$

\bigskip

\begin{cor} \label{Cor:3.2}  Let $U$ and $V$ be intervals contained in
	$\{ 1, 2, \ldots, \ell\}$.\\
(a) If $|U\cap V| \geq 2$ then $x_Uy_V=0 = y_Ux_V$. \\
(b) Otherwise $x_Uy_V$  and $y_Ux_V$ are non-zero.
\end{cor}
\bigskip

{\it Proof }   Let $G= U\cap V$.\\
(a) \ Assume $|G|\geq 2$. With the notation as in  Lemma \ref{lem:3.1} we have
$x_U=\psi^*x_G$ and $y_V=y_G\sigma$, and since $x_Gy_G=0$ it follows that $e_Uf_V=0$. 

(b) Suppose $|G|\leq 1$, then  $e_U$ and $f_V$ have disjoint supports, and
then the product is non-zero. Similarly $f_Ue_V$ is non-zero. 
$\Box$

\bigskip

\begin{defi} \label{def:3.3} \normalfont
We define for $1\leq r \leq \ell - 1$ and $r+s=\ell$, 
$$\ve_r:= \left\{\begin{array}{ll}  c_rx_{[1,s]}y_{[s+1, \ell]}  &   \mbox{ $r$ \  odd }\cr
                                                   c_ry_{[1, r]} x_{[r+1, \ell]}  &  \mbox{ $r$ \  even}
\end{array}
\right.
$$
	where $c_r$ is the non-zero scalar as described in 2.3 such that
	$\ve_r$ is an idempotent. 
\end{defi}

With these, we will show the following:

\begin{prop}\label{prop:3.4} We have $\ve_r\ve_u=0$ for $r\neq u$. Moreover
	$$\dim \ve_rH\ve_u = \left\{\begin{array}{ll} 0 & |r-u|\geq 2 \cr
	1 & u=r \pm 1 \cr
	2 & r=u.
\end{array}\right.
$$
\end{prop}

To proceed further, let 
$\lambda= (r, s)$ and $\mu = (u, v)$  be  compositions of $\ell$ with two non-zero parts. 
Slightly more general, we write 
$\ve_{\lambda}$ for one of $y_{[1, r]}x_{[r+1, \ell]}$ or
$x_{[1, r]}y_{[r+1, \ell]}$, and  we write $\ve_{\mu}$ for
one of $y_{[1,u]}x_{[u+1, \ell]}$ or $x_{[1,u]}y_{[u+1, \ell]}$.

 \begin{lemma} \label{lem:3.5} The space $\ve_{\lambda}H\ve_{\mu}$ is the span of the set
	 $\{\ve_{\lambda}T_{d}\ve_{\mu} \mid d \in \cD_{\lambda, \mu}\}$.
  \end{lemma}

{\it Proof}
 The space $\ve_{\lambda}H\ve_{\mu}$ is spanned by all $\ve_{\lambda}T_w\ve_{\mu}$ for $w\in W$.
Any such 
 $w$ can be written as  $w=h_1d h_2$ where $h_1\in W_{\lambda}$ and 
  and $h_2\in W_{\mu}$, and $d$ is a minimal length representative for the double coset containing $w$. Then $T_w = T_{h_1}T_{d}T_{h_2}$, and the claim follows.
     $\Box$

\parindent0pt

\begin{defi}\label{def:3.6}\normalfont Assume Proposition \ref{prop:3.4} and Lemma \ref{lem:3.5}. Let $1\leq r \leq \ell-2$. Define
	$$\alpha_r:= \ve_rT_d\ve_{r+1}, \ \ \beta_r = \ve_{r+1}T_{d^{-1}}\ve_r$$
where $d$ is the  distinguished double coset representative
	such that $\ve_rT_d\ve_{r+1}\neq 0$ (see \ref{prop:3.4} and \ref{lem:3.5}). 
\end{defi}

We will write down $\alpha_r, \beta_r$ explicitly below, and we will prove 
as the main result:

\bigskip

\begin{thm}\label{thm:3.7}  Let $B$ be the subalgebra of $H$ 
	generated by the set 
	$$\{ \ve_r\mid 1\leq r\leq \ell-1\}\cup
	\{ \alpha_r, \beta_r \mid 1\leq r\leq \ell-2\}.$$ 
	Then $B=\ve H\ve$ where $\ve = \sum_{r=1}^{\ell-1}\ve_i$ and  
	$B$ is a basic algebra for the principal block of $H$. \ 
After possibly rescaling some $\beta_r$ the elements  $\alpha_r, \beta_r$  satisfy the
	relations (3.1).
\end{thm}

 \vspace*{0.3cm}

\subsection{Towards $\ve_{\lambda}H\ve_{\mu}$.}

Let $\lambda = (r, s)$ and $\mu = (u, v)$ be compositions of $\ell$ with
non-zero parts.
We fix a distinguished  representative
$d\in \cD_{\lambda, \mu}$ and analyze the element $\ve_{\lambda}T_d\ve_{\mu}$, which  lies
in $H_{\lambda}T_dH_{\mu} = (T_d)(T_d)^{-1}H_{\lambda}T_d\cap H_{\mu}$. 
	By \cite{DJ}, Theorem 2.7, we have that
$$(T_d)^{-1}H_{\lambda}T_d\cap H_{\mu} = H_{\nu},$$
where $W_{\nu}$     is the Young subgroup generated by the set of basic transpositions contained in $d^{-1}W_{\lambda}d\cap W_{\mu}$. 
The proof of 2.7 in \cite{DJ} shows that if
$v=d^{-1}ud\in d^{-1}W_{\lambda}d\cap W_{\mu}$ for $u\in W_{\lambda}$ then 
$T_v = (T_d)^{-1}T_uT_d$. 
We will now compute $W_{\nu}$ (with $\nu = \nu(d)$), and also the Young subgroup $U$ of $W_{\lambda}$
such that $d^{-1}Ud= W_{\nu}$.

\medskip

(1) \ 
The elements $d\in \cD_{\lambda, \mu}$
are in bijection with the $2\times 2$ matrices $A$ with entries in $\bZ_{\geq 0}$ having row sums $u, v$ and column sums $r, s$, see for example Theorem 1.3.10 in \cite{JK}.
Take such a matrix
$A= \left(\begin{matrix} t & u-t\cr r-t & x \end{matrix}\right), 
$
where
$x = t + (s-u) = t+(v-r) (\geq 0).
$
Such an $A$ corresponds to the permutation $d$ which is,
written in  table form, 
$$\begin{matrix} 1 &  \ldots & t & t+1& \ldots & r &r+1 & \ldots & u+(r-t) & u+(r-t)+1  \ldots  \cr
1 &  \ldots & t & u+1 & \ldots &u+(r-t)  & t+1 & \ldots & u & u+(r-t)+1  \ldots 
\end{matrix}.
\leqno{(*)}$$
Note that if $t=r$ or $t=u$, then $d$ is the identity.
We observe that the matrix $A$ as above corresponds to $d\in \cD_{\lambda, \mu}$ if and only if its
transpose corresponds to $d^{-1}\in \cD_{\mu, \lambda}$. 

\begin{expl}\label{expl:3.8}
\normalfont 
(a) \ Let  $\lambda = (r,s)$ and $\mu = (s-1, r+1)$, where $r$ is even,         
$$\ve_{\lambda} =  y_{[1,r]}x_{[r+1, \ell]} \ \mbox{and} \ 
\ve_{\mu}  = x_{[1, s-1]}y_{[s, \ell]} 
$$
They  are non-zero scalar multiples of $\ve_r$ and $\ve_{r+1}$ respectively.  
We take the matrix $A$ and the corresponding $d$ as 
$$A = \left(\begin{matrix} 0 & s-1\cr r & 1\end{matrix}\right), \ \ 
d = \left(\begin{matrix} 1 & 2 & \ldots & r & r+1 & \ldots & \ell-1 & \ell\cr
                                        s & s+1 & \ldots & \ell-1 & 1 & \ldots & s-1 & \ell\end{matrix}\right).
$$
We will see later that with this, $\ve_rT_d\ve_{r+1}$ and $\ve_{r+1}T_{d^{-1}}\ve_r$ are  non-zero.

(b) Let $\lambda = (r, s)$ and $\mu = (s+1, r-1)$, where $r$ is even, and 
$$\ve_{\lambda} = y_{[1,r]}x_{[r+1, \ell]} \ \mbox{ and } \  \ve_{\mu} =  x_{[1, s+1]}y_{[s+2, \ell]}.
$$
They are non-zero scalar multiples of $\ve_r$ and $\ve_{r-1}$ respectively.
We take  the matrix $A$ and the corresponding permutation $d'$ as
$$A = \left(\begin{matrix} 1 & s\cr r-1 & 0\end{matrix}\right), \ \ 
d' = \left(\begin{matrix} 1 & 2 & \ldots & r & r+1 & \ldots & \ell\cr
1 &s+2& \ldots                            & \ell & 2 & \ldots & s+1\end{matrix}\right).
$$
	We will see later that $\ve_rT_{d'}\ve_{r-1}$ and $\ve_{r-1}T_{(d')^{-1}}\ve_r$ are non-zero.
\end{expl}

\medskip

\begin{lemma}\label{lem:3.9} We have
\begin{align*}
	W_{\nu} = & W_{[1, t]}\times W_{[u+1, u+(r-t)]} \times W_{[t+1, u]}\times
W_{[u+(r-t)+1, \ell]},\cr
	U= & W_{[1, t]}\times W_{[t+1, r]}\times W_{[r+1, u+(r-t)]}\times W_{[u+(r-t)+1, \ell]}.
\end{align*}
The relevant factors are trivial when $t=0$ or $t=r$ or $t=u$ or $u+(r-t)=\ell$.
 \end{lemma}

{\it Proof} \ 
To compute   $d^{-1}s_id$, we have
an elementary observation:
Let $d$ be the permutation which maps $j\mapsto a_j$, for $1\leq j\leq n$. Then for $1\leq i < n$ we have
that $d^{-1}s_id$ is the transposition $(a_i \ a_{i+1})$                        
Namely, $d^{-1}s_id$ is a transposition, and one sees directly that it takes $a_i\mapsto a_{i+1}$.
We apply this to $d$ as above, then 
$$d^{-1}s_id = \left\{\begin{array}{ll} s_i & 1\leq i\leq t-1\cr
	s_{i+(u-t)} & t+1\leq i\leq r-1\cr
	s_{i-(r-t)} & r+1\leq i\leq u+(r-t)-1\cr
	s_i & u+(r-t)+1 \leq i\leq \ell-1
\end{array}\right.
$$
Moreover $d^{-1}s_id$ is not a basic transposition for $i=t, r, u+(r-t)$.
With this, the Lemma follows.
$\Box$

\bigskip

(2) Now we will compute  an idempotent
$\ve_{\nu(d)} \in H_{\nu}$ such that $\ve_{\lambda}T_d = (\zeta T_d)\ve_{\nu(d)}$.
Consider first $\ve_{\lambda} = y_{[1,r]}x_{[r+1, \ell]}$.
We factorize $y_{[1,r]}$ and $x_{[r+1, \ell]}$ 
as described in Lemma \ref{lem:3.1} (taking $W_{\lambda} = \bigcup_{g\in \cD}gU$)  and we get 
$$\ve_{\lambda} = \zeta y_{[1,t]}y_{[t+1, r]}x_{[r+1, u+(r-t)]}x_{[u+(r-t)+1, \ell]}.$$
Here we write $\zeta$ for the product of the elements called $\psi^*$, and
$\sigma^*$, it is a linear combination of $T_g$ for $g\in \cD$. Then  

\begin{lemma} \label{lem:3.10} 
	(a) If $\ve_{\lambda}=y_{[1,r]}x_{[r+1, \ell]}$, then 
	$\ve_{\lambda}T_d = (\zeta T_d)\ve_{\nu(d)}$ where
	$$\ve_{\nu(d)}:= y_{[1,t]}y_{[u+1, u+(r-t)]}x_{[t+1, u]}x_{[u+(r-t)+1, \ell]} \leqno{(\dagger)}$$
(b) For $\ve_{\lambda}=x_{[1,r]}y_{[r+1, \ell]}$ we get the same formula
with $x, y$ interchanged.
\end{lemma}

{\it Proof} (a) \ With the above, 
\begin{align*} \ve_{\lambda}T_d &= 
	\zeta T_d(T_d)^{-1}y_{[1,t]}y_{[t+1, r]}x_{[r+1, u+(r-t)]}x_{[u+(r-t)+1, \ell]}T_d\cr
	&= \zeta T_dy_{[1,t]}(y_{[u+1, u+(r-t)]}x_{[t+1, u]})x_{[u+(r-t)+1, \ell]}\cr  
	& = (\zeta T_d)\ve_{\nu(d)}
\end{align*}

Part (b) is analogous.
$\Box$

\bigskip

\begin{prop}\label{prop:3.11} With the above setting, we have
	$\ve_{\lambda}T_d\ve_{\mu} = 0$ if and only if  $\ve_{\nu(d)}\ve_{\mu}=0$.
\end{prop}

\bigskip

{\it Proof } We must show that if $\ve_{\nu(d)}\ve_\mu$ is not zero then $(\zeta T_d)\ve_{\nu(d)}\ve_{\mu}$ is non-zero. 

Suppose  $\ve_{\nu(d)}\ve_{\mu}\neq 0$, then it is a linear combination of $T_w$ with $w\in W_{\mu}$. 
We have
$W_{\lambda} = \bigcup_{g\in \cD} gU$ and then
$d^{-1}W_{\lambda}d = \bigcup_{g\in \cD} d^{-1}gd W_{\nu}$, where $\cD$ is
the set of distinguished coset representatives.
Suppose $gdw= g'dw'$ for $g, g' \in \cD$,  and $w, w' \in W_{\mu}$. Then 
it follows that $d^{-1}(g')^{-1}dd^{-1}gd$ is in $d^{-1}W_{\lambda}d\cap W_{\mu} = W_{\nu}$. The $d^{-1}gd$ are coset representatives and it
follows that $g=g'$, and then also  $w=w'$.

Each $T_v$ in the support of  $\ve_{\nu(d)}\ve_{\mu}$ gives rise to a scalar multiple of 
$T_gT_dT_w$ in $(\zeta T_d)\ve_{\nu(d)}\ve_{\mu}$, and we  have  
$T_gT_dT_w= T_{gdw}$, using that $d$ is a distinguished double coset representative. 
As we  have just proved, there is no repetition amongst  these elements $gdw$ 
as $g$ varies through  $\cD$  and $w$ varies 
through $W_{\mu}$, and hence we have
an expression of $(\zeta T_d)\ve_{\nu(d)}\ve_{\mu}$ in terms of
the basis of $H$, with non-zero coefficients.
$\Box$

\medskip

\begin{cor}\label{cor:3.12} With the same
        notation, the set $\{ T_{gdw}\mid w\in W_{\mu}, g\in \cD\}$ is linearly
        independent.
\end{cor}

This is part of the above proof.

\medskip

\subsubsection{Proof of Proposition \ref{prop:3.4}} 

Let $\lambda = (r, s)$ and $\mu = (u, v)$ with each of $r, s, u, v$ non-zero.

(I) \ Assume  that $r, u$ have the same parity. 
We consider first the case when 
$\ve_{\lambda} = y_{[1,r]}x_{[r+1, \ell]}$ and 
$\ve_{\mu} = y_{[1,u]}x_{[u+1,\ell]}$.
Let $d\in \cD_{\lambda, \mu}$, then by Lemma \ref{lem:3.10} we have  $\ve_{\lambda}T_d = (\zeta T_d)\ve_{\nu(d)}$ where
$\ve_{\nu(d)}$ is the element $(\dagger)$ of \ref{lem:3.10}.
Note  that any two factors of $\ve_{\nu(d)}$ commute.\\
(i)  \ Suppose $t+1<u$ then $x_{[t+1,u]}y_{[1, u]}=0$, and hence $\ve_{\nu(d)}\ve_{\mu}=0$. \\
(ii) \ Suppose $1<(r-t)$ then $y_{[u+1, u+(r-t)]}x_{[u+1, \ell]}=0$ and again $\ve_{\nu(d)}\ve_{\mu}=0$.\\
Recall that $t\leq u$, so this leaves $t+1=u$ or $t=u$, and $t+1=r$, or $t=r$. We assume that $u$ and $r$ have the same parity, therefore $u=r$ and it is
equal to $t$ or $t+1$.
In particular $\ve_{\lambda} = \ve_{\mu}$, and consequently we have $\ve_{\lambda}H\ve_{\mu}=0$ for $\lambda\neq \mu$ in this case.

\medskip

 If $r=t$ then $d$ is the identity  and $\ve_{\lambda}T_d\ve_{\lambda}$ is non-zero
(by 2.3). Suppose $r=t+1$, then  $\ve_{\nu(d)}\ve_{\lambda} = y_{[1, r-1]}x_{r+1,\ell]}y_{[1,r]}x_{[r+1,\ell]}$ which by 2.3 is a non-zero scalar multiple of  $\ve_{\lambda}$, and hence
$\ve_{\lambda}T_d\ve_{\lambda}$ is non-zero (by \ref{prop:3.11}). These are two non-zero elements supported on different double cosets, so they are
linearly independent.
Hence $\ve_{\lambda}H\ve_{\lambda}$ is 2-dimensional. 

\medskip

Now assume $\ve_{\lambda} = x_{[1,r]}y_{[r+1, \ell]}$ and $\ve_{\mu} = x_{[1,u]}y_{[u+1,\ell]}$. Then  the same proof works with $x, y$ interchanged, and we 
get again that $\ve_{\lambda}H\ve_{\mu}$ is zero for $\lambda\neq \mu$ and is 2-dimensional otherwise.
 This proves Proposition \ref{prop:3.4} when $r, u$ have the same parity.

\bigskip

(II)  \ Now we deal with idempotents  $\ve_r$ for $r$ even and $\ve_v$ and
$v$ odd. 
Suppose  $\ve_{\lambda} = y_{[1,r]}x_{[r+1, \ell]}$ and $\ve_{\mu} = x_{[1,u]}y_{[u+1, \ell]}$.  Let $d\in \cD_{\lambda, \mu}$, then 
$\ve_{\lambda}T_d = (\zeta T_d)\ve_{\nu(d)}$ with $\ve_{\nu(d)}$ as in $(\dagger)$ of \ref{lem:3.10}. \\
(i) \ 
Suppose $t\geq 2$,  then  $y_{[1, t]}x_{[1,u]}=0$ and hence $\ve_{\nu(d)}\ve_{\mu}=0$. \\
(ii) \ Suppose  $u+(r-t)+1 <  \ell$, then  $x_{[u+(r-t)+1, \ell]}y_{[u+1, \ell]} = 0$. and hence 
$\ve_{\nu(d)}\ve_{\mu}=0$. \\
Since $t\leq u$, this  leaves $t=0$ or $t=1$, and $u+(r-t) \geq \ell-1$. 
Recall from the first part of Section 3.2, 
the bottom left entry of the matrix $A$ is   
$x= t+(s-u) = t+(v-r)\geq 0$. \\

(a) Assume $t=0$,  then $s\geq u$ and $v\geq r$ and in fact $v>r$ since $v, r$ have
different parity (and then $s>u$).
In this case we have $u+r =  \ell-1$ since  $r\neq v$.
It follows that $r+1=v$ (hence $\ve_v= \ve_{r+1}$). 
Moreover 
$$\ve_{\nu(d)}\ve_{\mu} = y_{[u+1, \ell-1]}x_{[1,u]}\ve_{\mu}
\leqno{(\ve_v=\ve_{r+1})}$$
which is a non-zero scalar multiple of $\ve_{\mu}$. 

(b) Assume  $t=1$ then similarly $u+(r-1)=\ell$ and $r-1=v$
(hence $\ve_v=\ve_{r-1}$). Moreover
$$\ve_{\nu(d)}\ve_{\mu} = y_{[u+1, \ell]}x_{[2, u]}\ve_{\mu}
\leqno{(\ve_v=\ve_{r-1})}$$ 
which again is a non-zero scalar multiple of $\ve_{\mu}$.
Hence in this case $\ve_rH\ve_v$ is 1-dimensional if $v=r\pm 1$ and is zero otherwise. 
\medskip

(II') \ Consider now idempotents $\ve_s$ and $\ve_u$ for $s$ odd and $u$ even; we  take 
 $\ve_{\lambda} = x_{[1,r]}y_{[r+1, \ell]}$ and $\ve_{\mu} = y_{[1,u]}x_{[u+1, \ell]}$, and let $d\in \cD_{\lambda, \mu}$.
With the same reasoning as in (II) we get $\ve_{\nu(d)}\ve_{\mu}\neq 0$ only  
 for $t=0$ and $t=1$.\\ 
 (a) When $t=0$, it follows that  $\ve_u = \ve_{s-1}$. Moreover
 $$\ve_{\nu(d)}\ve_{\mu} = x_{[u+1, \ell-1]}y_{[1,u]}\ve_{\mu} \leqno{(\ve_u=\ve_{s-1})}
 $$
 which is a non-zero scalar multiple of $\ve_{\mu}$.
 
 (b) When $t=1$, we get $\ve_u=\ve_{s+1}$. Moreover 
$$\ve_{\nu(d)}\ve_{\mu} = x_{[u+1, \ell]}y_{[2,u]}\ve_{\mu} \leqno{(\ve_u=\ve_{s+1})}
 $$
and this also is a non-zero scalar multiple of $\ve_{\mu}$. 
Hence  $\ve_sH\ve_u$ is 1-dimensional for 
$\ve_u= \ve_{s\pm 1}$ and is zero otherwise.
This completes the proof of Proposition 3.4.
$\Box$

 \bigskip

 \begin{rem}\label{rem:3.13}\normalfont We collect information from the
	 above proof.\\
	(1) The matrix $A$ in
	 (II)(a) corresponds to  $d$ in Example \ref{expl:3.8}(a), this
	 gives $\alpha_r$ as in Definition \ref{def:3.6} for $r$ even.\\
	 (2)
	 The matrix $A$  in 
	  (II)(b)  corresponds to  $d'$ in Example \ref{expl:3.8}(b), this
	 gives $\beta_{r-1}$ as in Definition \ref{def:3.6} for $r$ even.\\
	 (3) The  permutation associated to the matrix in part (II')(a) 
	 (the inverse of the permutation as in (1) with $s-1$ instead of $r$)  defines $\beta_{s-1}$.\\
	 (4) The  permutation associated to the matrix in part (II')(b) 
           (the inverse of the permutation as in (2) with $s=r-1$)  defines $\alpha_s$.
 \end{rem}

\subsection{Arrows and relations}

\begin{lemma}\label{lem:3.14} The algebra $\ve_1H\ve_1$ has basis $\ve_1, x_{[1, \ell]}$. It is local, isomorphic to $K[X]/(X^2)$.
\end{lemma}

\medskip

{\it Proof} \ 
The element $x_{[1, \ell]}$ spans the trivial $H$-module, and
its square is zero.
One checks that
$x_{[1, \ell]} = \ve_1x_{[1, \ell]}\ve_1$. 
Clearly $x_{[1,\ell]}$ and $\ve_1$ are linearly independent.
We have proved in Proposition \ref{prop:3.4} 
that $\ve_1H\ve_1$ is 2-dimensional, and it follows that
it has basis $\{ \ve_1, x_{[1,\ell]} \}$.
The Lemma follows.
$\Box$

\bigskip

We will now prove  Theorem \ref{thm:3.7}.    Let $\alpha_r, \beta_r$ be as in 
Definition \ref{def:3.6}. 

\medskip
(1) \ Clearly $\alpha_r\alpha_{r+1}=0$ since it is an element
in $\ve_rH\ve_{r+2}$ which is zero by Proposition \ref{prop:3.4}, similarly
$\beta_{r+1}\beta_r=0$. 

\medskip
(2) \  The products $\alpha_r\beta_r$ and $\beta_r\alpha_r$ are
	non-zero:\\
We  fix $r$ and write $\alpha = \alpha_r$ and $\beta=\beta_r$.
To prove $\alpha\beta \neq 0$, it suffices to 
show  $f(\alpha\beta, 1) \neq 0$ 
where $f$ is the symmetrizing
form as in 2.4. 
Using the symmetry property, and the fact that the $\ve$ are idempotent, 
and also that $\ve_rT_d\ve_{r+1} = c (\zeta T_d)\ve_{r+1}$ where $c$ is
a non-zero scalar, we get 
\begin{align*}
	f(\alpha\beta, 1) = & f(\ve_rT_d\ve_{r+1}T_{d^{-1}}, \ve_r) = f(\ve_r^2T_d\ve_{r+1}T_{d^{-1}}, 1) \cr
	=& f(\ve_rT_d \ve_{r+1}, T_{d^{-1}}) = 
	 (\ve_rT_{d}\ve_{r+1},  T_d)
\end{align*}
where $(-,-)$ is the bilinear form as in 2.4.  We have that 
$\ve_rT_d\ve_{r+1} = (\zeta T_d)\ve_{\nu(d)}\ve_{\mu}$ with the notation
as in Proposition \ref{prop:3.11} (with $\ve_{r+1}$ a non-zero scalar multiple of $\ve_{\mu}$).
We know that $(\zeta T_d)\ve_{\nu(d)}\ve_{\mu}$ is nonzero, so by Corollary \ref{cor:3.12}, it
is a unique linear combination of the set $\{ T_{gdw}\mid g\in \cD, w\in W_{\mu}\}$. Therefore 
$$((\zeta T_d)\ve_{\nu(d)}\ve_{\mu}, T_d) = cq^{l(d)}
$$
where $c$ is the coefficient of the identity in $\ve_{\nu(d)}\ve_{\mu}$ and this is non-zero (by 2.2).

\medskip

(3) \ $\alpha_r\beta_r\alpha_r=0$ and $\beta_r\alpha_r\beta_r=0$, in particular each $\ve_rH\ve_r$ contains a non-zero nilpotent element, and hence
is a 2-dimensional local algebra:\\
To prove this, we
start with $r=1$ and use Lemma \ref{lem:3.14}. It suffices to show that
$\alpha_1\beta_1$ is not a unit  in $\ve_1H\ve_1$.
We have 
$$x_{[1, \ell]}\alpha_1 = x_{[1, \ell]}\ve_1T_d\ve_2 = x_{[1, \ell]}T_d\ve_2
= q^{l(d)}x_{[1, \ell]}\ve_2 =0, 
$$
hence $x_{[1, \ell]}\alpha_1\beta_1=0$. But $x_{[1, \ell]}\neq 0$ and therefore $\alpha_1\beta_1$ is not a unit. 
We deduce that  $\alpha_1\beta_1$ must be nilpotent and is therefore  a scalar multiple of $x_{[1, \ell]}$. By  the above calculation, 
we see directly that $\alpha_1\beta_1\alpha_1=0$.
We can also deduce $\beta_1\alpha_1\beta_1=0$: \ It is of the form $a \beta_1$ for $a\in K$ and then $0 = (\alpha_1\beta_1)^2 = a(\alpha_1\beta_1)$
and $a=0$.  This  shows that 
$\beta_1\alpha_1 \in \ve_2H\ve_2$ is nilpotent and then the algebra
$\ve_2H\ve_2$ must be local as well. 
By induction on $r$, repeating the arguments, the claim follows.
\medskip

(4) 
 \ Up to rescaling, the commutativity relations hold: We use that by the previous part, the radical of $\ve_rH\ve_r$ is 1-dimensional. 
 The non-zero elements $\beta_r\alpha_r$ and $\alpha_{r+1}\beta_{r+1}$ 
are both in the 1-dimensional radical of  $\ve_rH\ve_r$ and hence must be scalar
multiples of each other. 
So we can take  $\beta_2'= a_2\beta_2$ for $0\neq a_2\in K$ so that
$\beta_1\alpha_1 = \alpha_2\beta_2'$.

Inductively,  suppose we have scaled arrows so that for $i\leq r$ 
the commutativity relation 
$\beta_{i-1}'\alpha_{i-1}= \beta_i'\alpha_i$ holds. If $r < \ell-2$ then 
$\beta_{r+1}'=a_{r+1}\beta_{r+1}$ so that
$\beta_r'\alpha_r = \alpha_{r+1}\beta_{r+1}'$.

\medskip

(5) \ The relations (1) to (4) determine the algebra: We can write $B$ as
 $\ve H \ve$ where $\ve = \sum_{r=1}^{\ell-1} \ve_r$.
Our computations 
show that it has dimension $2(\ell-1) + 2(\ell-2) = 4(\ell-1)-2$. This is also
the dimension of the algebra $KQ/I$ in (3.1).
One defines an algebra map $KQ \to B= \ve H\ve$ by mapping the elements
$e_r\in KQ$ (corresponding to paths of length zero) to $\ve_r$, and taking
arrows $\alpha_r, \beta_r$  to the elements 
in $\ve H\ve$ we called $\alpha_r$ and the rescaled $\beta_r$ , and extend
to products and linear combinations.
It is surjective, 
and the ideal $I$ is in the kernel. Then by dimensions, $I$ is the kernel
and we have shown that $\ve H\ve$ is the algebra as defined above.
$\Box$

\bigskip

Recall that the socle of a module is the largest semisimple submodule, and the top is the
largest semisimple quotient.
\medskip

\begin{cor}\label{cor:3.15}  (a) Each   $\ve_rH$ is indecomposable projective.
	For $2\leq r \leq \ell-2$ Its radical has basis $\{ \alpha_r, \beta_{r-1}, \alpha_r\beta_r\}$. The radical of $\ve_1H$ has basis $\{\alpha_1, \alpha_1\beta_1\}$
	and the radical of $\ve_{\ell-1}H$ has basis $\{\beta_{\ell-2}, \beta_{\ell-2}\alpha_{\ell-2}\}$.\\ 
	(b) The socle of $\ve_rH$ is simple, spanned by $\alpha_r\beta_r$ (or $\beta_{\ell-2}\alpha_{\ell-2}$).\\
(c) The simple $\ve H\ve$-modules are precisely the 1-dimensional quotients $\ve_rH/ {\rm rad} \ve_rH$.
The socle of $\ve H \ve$ is spanned by the $\alpha_r\beta_r$ ($1\leq r\leq \ell-2$) together with $\beta_{\ell-2}\alpha_{\ell-2}$. 
\end{cor}

\bigskip

{\it Proof }
Since  $\ve_r$ is an idempotent, the module $\ve_rH$ 
is projective.
Its endomorphism algebra is
$\ve_rH\ve_r$, and we have proved that it is local. Hence
$\ve_rH$ is indecomposable.
All other statements follow from the proof of Theorem \ref{thm:3.7}
$\Box$

\begin{rem}\label{rem:3.16} \normalfont  (1) The module $\ve_rH$ is 
the projective cover
	of $D^{ ((\ell-r)+1, 1^{r-1})}$ for $r\geq 1$. This follows from the decomposition number approach to such blocks,
	see for example \cite{Ge}, together with the observation that $\ve_1H$ is the projective cover of the
	trivial module.
\end{rem}

\bigskip

\begin{expl}\label{expl:3.17} \normalfont Assume $\ell=3$.
\  The algebra $H=H_q(3)$ has dimension $6$ and hence is equal to 
	its  basic algebra $B$ by dimensions.
The presentation as a Hecke algebra is  not compatible with
the presentation (3.1). The generators in our presentation are
$\ve_1 = c_1x_{[1,2]}$ and $\ve_2 = c_2y_{[1, 2]}$
(with $c_1, c_2$ non-zero scalars), and the arrows are
$\beta_1=\ve_2T_2\ve_1$ and $\alpha_1 = \ve_1T_2\ve_2.
$
\end{expl}

\bigskip

\section{Homological properties}

\subsection{Controlling projectivity}\label{subsec:4.1}

Given a module, a basic question  is to determine whether or not it
is projective. 
More generally, an indecomposable module of a symmetric, or even selfinjective, algebra is
projective if and only if it the socle of the algebra does not annihilate
the module. For an algebra which is not basic, this is not so easy to
verify.
In our setting, we can use this
and give an easy criterion.

We show  (with  $\ell \geq 3$)  that there is an associated truncated polynomial algebra which directly detects projectivity.

 Consider the basic algebra $A$ for the principal block of $H_q(\ell)$ as in Theorem \ref{thm:3.7}.
We choose and fix non-zero elements $c_r\in K$ for $1\leq i\leq \ell - 2$ such that $c_r+ c_{r+1}\neq 0$. (Such elements exist, $K$ has at least two non-zero elements.)

\begin{lemma} \label{lem:4.1} Let
$$\wt{\alpha}:= \sum_{r=1}^{\ell - 2}\alpha_r\ \mbox{ and }
	\wt{\beta}:= \sum_{r=1}^{\ell - 2} c_r\beta_r,  \mbox{ and }
\sigma:= \wt{\alpha} + \wt{\beta}.$$
Then  $\sigma^2$ generates the socle of $A$ as an $A$-module.
	Furthermore  $\sigma^3=0$.
\end{lemma}

\medskip

{\it Proof } First observe that $\wt{\alpha}^2=0$ and $\wt{\beta}^2=0$. Then
$$\sigma^2
= c_1\alpha_1\beta_1 + \sum_{r=2}^{\ell-3} (c_r+c_{r+1})\alpha_r\beta_r + c_{\ell - 2}\beta_{\ell-2}\alpha_{\ell-2}.
$$
All scalar coefficients are non-zero by construction. 
Moreover $\ve_r\sigma^2\ve_r$ spans the socle of $\ve_rH$, hence
$\sigma^2$  spans the socle of $A$ and  $\sigma^3=0$.
$\Box$

\bigskip

Let $R$ be the subalgebra of $A$ generated by $\sigma$, this is isomorphic to
$K[X]/(X^3)$. 
Now take a block of an $\ell$-parabolic subalgebra which is not simple.
This is Morita equivalent to the tensor product $A^{\otimes m}$ of copies of $A$ for some $m\geq 1$, let $B$ be this
algebra.  Let $\sigma_i\in B$ be the tensor  where at the $i$-th place take $\sigma$, and all other factors
are the identity. Then let $R_B$ be the subalgebra of $B$ generated by $\sigma_1, \cdots, \sigma_m$. The $\sigma_i$ commute, and
$\sigma_i^3=0$ so that  $R_B$ is a truncated polynomial algebra.

\medskip

\begin{prop}\label{prop:4.2}  Assume $M$ is an indecomposable $B$-module. Then $M$ is projective
	as a $B$-module if and only if
$(\sigma_1\sigma_2\cdots \sigma_m)^2$ does not
annihilate $M$.
\end{prop}

{\it Proof} \ 
(1) Suppose $M$ is projective, then $M\cong eB$ where $e=\ve_{i_1}\otimes \ldots \ve_{i_m}$ and each $\ve_{i_j}$ is one of the idempotents
$\ve_1, \ldots \ve_{\ell-1}$ in $A$. 
We have 
$$
	e(\sigma_1\sigma_2\cdots \sigma_m)^2 = \ve_{i_1}\sigma^2\otimes \ldots \otimes \ve_{i_m}\sigma^2.
$$
This spans the socle of $M$ since 
$\ve_j\sigma^2$ spans the socle of $\ve_jA$ for each $j$. In particular
$M(\sigma_1\cdots \sigma_m)^2$ is non-zero.

(2) If $M$ is indecomposable and $M(\sigma_1\cdots \sigma_m)^2\neq 0$ then there
is $m = me\in M$ and an idempotent $e=\ve_{i_1}\otimes \ldots \ve_{i_m}$ as
above with $me(\sigma_1\cdots \sigma_m)^2 \neq 0$. Then one shows that
the submodule $mB=meB$ of $M$ is isomorphic to $eB$, using that the socle of
$eB$ is simple. Now $M$ is indecomposable  and $eB$ is projective but also is injective, hence $M \cong eB$ and
$M$ is projective.
$\Box$

\bigskip

We would like to know when   the algebra $B$ is projective as a module over $R$.

\medskip

\begin{lemma}\label{lem:4.3} Let $A$ the basic algebra of the principal block of $H_q(\ell)$ and $B=A^{\otimes m}$ for $m\geq 1$.
Then $B$ is projective as an $R_B$-module if and only if $\ell = 3$.
\end{lemma}

{\it Proof }
(1) \ We reduce to $m=1$: Suppose $A$ is not projective as an $R$-module,
say it is $P\oplus M$ with $P$ projective, and where $M$ is non-zero and has no
projective summand. Then $M\sigma^2=0$, and the summand $M^{\otimes m}$ of $B$ as an $R$-module is not projective. Hence  $B$ is not projective as a module for $R_B$. 

(2) Now consider $m=1$.
Assume first that $\ell=3$, we can see directly that $A$ is projective as an 
$R$-module by checking  that $A= \ve_1R \oplus \ve_2R$.

Now assume  $A$ is projective as an $R$-module. We have
$A=\oplus_{i=1}^{\ell-1} \ve_iA$,  and each summand is  $R$-invariant.
It follows that each $\ve_iA$ must be projective as an $R$-module, and
hence must have dimension divisible by $3$.
Suppose $\ell>3$, then $\ve_2A$ has dimension $4$, a contradiction.
Therefore $\ell = 3$. 
$\Box$

\bigskip

This suggests that $R$ may not be so useful in general.

\bigskip

\subsection{Support varieties}

From now we assume that $K$ is algebraically closed.
To control projectivity of modules, and understand  some large-scale
behaviour, one can often exploit cohomological support varieties.
Hecke algebras are not Hopf algebras but one can 
use Hochschild cohomology. We recall relevant definitions and facts, 
here $\La$ is a finite-dimensional algebra.
We refer to \cite{EHSST} for details.

\medskip

\begin{defi} The algebra  $\La$ satisfies the finite generation hypothesis
	(Fg) if the Hochschild cohomology $HH^*(\La)$ is Noetherian and
	 moreover ${\rm Ext}^*_{\La}(\La/{\rm rad} \La, \La/{\rm rad}\La)$
	is a finitely generated $HH^*(\La)$-module. 
\end{defi}

With this, any $\La$-module $M$ has a support variety $V(M)$ whose properties are similar to that of support varieties for group representations. 
The dimension of $V(M)$ is equal to the complexity of $M$, in particular 
$V(M)$ is trivial if and only if $M$ is projective. Furthermore,
it gives information on Auslander-Reiten components. 

\begin{defi} \normalfont Assume $\cX$ is a subcategory of  mod-$\La$ which is
	the union of stable Auslander-Reiten components 
	(up to projectives). We say
	that there are enough periodic modules for $\cX$ if for every
	non-projective indecomposable module $M$ in $\cX$ there is 
	some $\Omega$-periodic $\La$-module $W$ such that
	$\ul{\rm Hom}_{\La}(W, M)\neq 0$.  If so, then 
the tree class of any component is one of the following 
(see for example \cite{Be}, or \cite{BEH} 2.11 and 2.12)
$$\mbox{{\it  
 Dynkin of type ADE, or  Euclidean, or an infinite tree}}
\ A_{\infty}, \  D_{\infty}, \  A_{\infty}^{\infty}.   \leqno{(\cT)}
$$
\end{defi}

\bigskip

We consider the case when $\La = H_{\rho}$, an $\ell$-parabolic Hecke
algebra. Then we have:

\medskip

\begin{lemma} \label{lem:4.6} The algebras $H_{\rho}$ satisfy (Fg). In particular\\
(a) $H_{\rho}$-modules have finite complexity.\\
(b) The category of $H_{\rho}-$modules has enough periodic modules. \\
	(c) The tree class of any  stable AR component belongs to $(\cT)$.
\end{lemma}

\bigskip
{\it Proof } 
The condition (Fg) is Morita invariant, so it suffices to show 
that the basic algebra for each block of $H_{\rho}$ satisfies
(Fg). For the basic algebra $A$ of $H_q(\ell)$  this follows from  \cite{ES} (the algebra has radical cube zero and is Koszul, one checks that its
Ext algebra is finitely generated over its graded center). Then by 
\cite{BO} (Corollary 4.8), every tensor power of $A$
satsifies (Fg). 
Parts (b) and (c) follow directly, by applying results from \cite{EHSST}.
$\Box$

\begin{rem}\normalfont (1) \ So far this deals with $H_{\rho}$. 
        If the field has
  characteristic zero, Linckelmann's work \cite{Li} implies that
	even $H$ satisfies 
  (Fg).
It is open whether it is true for non-zero characteristic, and
there are many open questions. For example it is not even known whether
	the trivial module for $H_q(6)$ when $\ell=3$ and $p=2$  has finite complexity.

(2) \ Over characteristic zero, an explicit presentation
	of the cohomology of $H$, that is, of
	${\rm Ext}^*_H(K, K)$ with $K$ the trivial module,
	was determined in \cite{BEM}. This is used by \cite{NX} to develop
	a support variety theory.
\end{rem}

\medskip

\subsection{Rank varieties} \label{subsec:4.3}
For some algebras, modules have rank varieties (and often it is known that they are essentially the same
as the support varieties). Rank varieties were first introduced
for group algebras of elementary abelian $p-$groups over fields of characteristic $p$ by J. Carlson
\cite{C}. This was generalized to quantum complete intersections  in \cite{BEH}.

These include truncated polynomial algebras, in particular
the algebra $R$ we have
constructed earlier, and $H_{\rho}$ itself when $\ell=2$. 
Rank varieties can be used directly to
construct enough periodic modules, see \cite{BEH}; see also 
\cite{Sch}. 
A recent result in  \cite{AvIy} 
introduces a different type of rank variety.

\medskip

In \cite{Sch},  S. Schmider shows, using techniques similar to \cite{BEH}:
\medskip

\begin{thm} Suppose $\ell \geq 3$ and  $p$ does not divide $\ell-1$. 
Then all components in a block of wild representation type
of $H_{\rho}$ have tree class $A_{\infty}$.
\end{thm}

\medskip
The condition comes because he uses skew group rings.

\subsection{Beyond $\ell$-parabolic subalgebras}

Let $H_{\rho}$ be a maximal   $\ell$-parabolic subalgebra.
Define  $\cX$ to be the full subcategory of mod-$H$ with objects the $H$-modules
$M$ whose restriction to $H_{\rho}$ is not projective.

\begin{prop} (a) \ The category $\cX$ is the union
	of stable Auslander-Reiten components. \\
(b) \ There are enough periodic modules for $\cX$, hence each tree class belongs
	to the list $(\cT)$.
\end{prop}

\bigskip
{\it Proof} \ (a) Clearly $\cX$ is  closed under syzygies, and under
Auslander-Reiten translation $\tau (\cong  \Omega^2$).
Let $0\to \tau(M)\to E\to M\to 0$ be an almost split sequence, and assume
$M$ is in $\cX$, then also $\tau(M)$ is in $\cX$. Take an indecomposable
non-projective summand $E'$ of $E$,
then by general theory there is an almost split sequence
$$0\to \tau(E')\to U\oplus \tau(M) \to E'\to 0\leqno{(*)}$$
(for some module $U$). Assume for a contradiction $E'\not\in \cX$, that is $E'_{H_{\rho}}$ is projective, then the restriction of (*) to $H_{\rho}$ is split. 
As well $\tau(E')$ is not in $\cX$ and so $E'\oplus \tau(E')$ is projective
restricted to $H_{\rho}$, and we have
the contradiction that $\tau(M)$ is  projective as a module for $H_{\rho}$.

\ \ (b) Let $M$ be indecomposable and not projective, and $M\in \cX$, then there is a periodic $H_{\rho}$-module
$W$ such that $\ul{\rm Hom}_{H_{\rho}}(W, M)\neq 0$ (by 4.6).
Now, $H$ is projective as a module for $H_{\rho}$, so 
$\ul{\rm Hom}_H(W\otimes_{H_{\rho}}H, M)\neq 0$. 
Furthermore, the  module $W\otimes_{H_{\rho}}H$ is periodic up
to projective summands. 
$\Box$

\bigskip

In general, it is not known which tree classes for $H$-modules occur,
except for blocks of finite type (with tree class $A_{\ell-1}$), and for
tame type (with one Euclidean component and otherwise tubes, with tree class
$A_{\infty}$, see \cite{EN}.)
\medskip

\begin{expl} \normalfont 
	Let $H= H_q(2\ell)$ for $\ell \geq 3$, and assume char$(K)=2$. Let  
$\rho = (\ell, \ell)$, then  $H_{\rho}$ is maximal $\ell$-parabolic.
We will  show that there
is  a module $V$ which is not projective but the restriction of
$V$ to $H_{\rho}$ is projective. 
We identify   $V$ 
as a subquotient of the module $x_{\rho}H$, the $q$-permutation
module $M^{\rho}$
(in the terminology of \cite{DJ}).

\medskip

Using 2.5, 2.6 and 3.3 from \cite{DJ}, noting that
	$K\cong M^{(2\ell)}$, one can see that 
(up to scalar multiples) there is a unique non-zero homomorphism
$\vf: K\to M$, and also an essentially unique non-zero homomorphism
$\psi: M\to K$.
Furthermore,  these split
when restricted to $H_{\rho}$.

The composition $\psi\circ \vf : K\to K$ must be zero, since
	otherwise $K$ would be  a direct summand of $M^{\rho}$. But this is not the case;  
	the vertex of $Y^{(n)}$ is always a  maximal $\ell-p$ parabolic
	subalgebra (see \cite{DD}), which in our case is $H$.
So let $U_1:= {\rm Im}(\vf)$ and $U_2:= {\rm Ker}(\psi)$, then we have
the chain of submodules
$0\subset U_1\subset U_2\subset M^{\rho}$
with $U_1$ and $M^{\rho}/U_2$ both isomorphic to $K$.
\medskip

Let $V:= U_2/U_1$. One can show that $V$ is not projective, by analysing the component
	in the principal block.  By  the Mackey decomposition (see \cite{DJ} 2.7)
	the restriction of $M^{\rho}$ to $H_{\rho}$ is isomorphic to $K \oplus K \oplus P$
	with $P$ projective. 
Since  the maps $\vf, \psi$ split on restriction to $H_{\rho}$, it
follows that $V$ as a module for $H_{\rho}$ is isomorphic to the
projective module $P$.
\end{expl}

\end{document}